\documentclass[12pt]{article}
\font\smallit=cmti12
\usepackage{theorem,amssymb,amsmath}
\theorembodyfont{\rmfamily}

\newtheorem{thm}{Theorem}[section]
\newtheorem{lemma}[thm]{Lemma}
\newtheorem{cor}[thm]{Corollary}
\newtheorem{prop}[thm]{Proposition}

\makeatletter

\renewcommand\section{\@startsection {section}{1}{\z@}%
{-30pt \@plus -1ex \@minus -.2ex}%
{2.3ex \@plus.2ex}%
{\normalfont\normalsize\bfseries}}

\renewcommand\subsection{\@startsection{subsection}{2}{\z@}%
{-3.25ex\@plus -1ex \@minus -.2ex}%
{1.5ex \@plus .2ex}%
{\normalfont\normalsize\bfseries}}
\renewcommand{\@seccntformat}[1]{\csname the#1\endcsname. } 

\makeatother

\begin{document}

\date{}
\title{\bf \normalsize \uppercase{Bounds on Van der Waerden Numbers \\and Some
Related Functions}}
\author {\bf \normalsize Tom Brown \\
{\smallit Department of Mathematics, Simon Fraser
University, Burnaby, BC V5A 1S6} 
\and
\bf \normalsize Bruce M. Landman\\
{\smallit Department of Mathematics, University of West Georgia,
Carrollton, GA 30118}
\and 
\bf \normalsize Aaron Robertson \\
{\smallit Department of Mathematics,
Colgate University, Hamilton, NY 13346}
}
\maketitle

\begin{abstract}
\noindent
For positive integers $s$ and
$k_1, k_2, \dots, k_s$, let $w(k_1,k_2,\dots,k_s)$ be the minimum integer $n$
such that any
$s$-coloring $\{1,2,\dots,n\} \rightarrow \{1,2,\dots,s\}$ admits  a $k_i$-term
arithmetic progression of color $i$
for some $i$, $1 \leq i \leq s$.  In the case when
$k_1=k_2=\cdots=k_s=k$ we simply write $w(k;s)$.
 That such a minimum integer exists
follows from van der Waerden's theorem on arithmetic progressions.
	In the present paper we give a lower bound for $w(k,m)$
for each fixed  $m$.  We include a table with values of  $w(k,3)$
which match this lower bound closely for $5 \leq k \leq 16$.  We also give an upper
bound for $w(k,4)$, an upper bound for  $w(4;s)$, and a lower bound for  $w(k;s)$ for
an arbitrary fixed  $k$.  We discuss a number of other functions that are closely
related to the van der Waerden function.

\end{abstract}

\baselineskip=14pt

\section{Introduction}

Two fundamental theorems in combinatorics are van der Waerden's
Theorem [19] and Ramsey's Theorem [15].
The theorem of van der Waerden  says, in particular, that for any two
given positive integers $k$ and $m$, there
exists a least positive integer $n = w(k, m)$ such that whenever
the integers in $[1, n] = \{1, 2, \dots , n\}$ are colored with two
colors (i.e., partitioned into two sets), there is either a $k$-term
arithmetic progression of the first
color (i.e., contained in the first set) or an $m$-term arithmetic
progression of the second color (i.e., contained in the second set).

Similarly, Ramsey's Theorem has an associated ``threshold" function
$R(k, m)$ (which we will not define here).  This function  satisfies the inequality
\[ R(k, m)  \leq  R(k-1, m)  +  R(k, m-1),\]
which leads to an upper bound on $R(k, m)$ that is not so much larger than
the best known lower bounds on $R(k, m)$ obtained by other means. 
Furthermore, the order of magnitude of $R(k, 3)$ is known to be $\frac{k^2}{\log
k}$ [10].

    For van der Waerden's function $w(k,m)$ there is no corresponding recursive
inequality known,
    and the order of magnitude of $w(k,3)$ is not known.  The best known lower and
upper bounds on
$w(k, k)$
     are
\[ (k-1)2^{(k-1)}  \leq w(k, k) < 2^{2{^{2^{2^{2^{(k+9)}}}}}}, \]
the lower bound known only when $k-1$ is prime.  The lower bound is
due to Berlekamp [2] and the upper bound to Gowers [6]. Narrowing
this gap is a fundamental problem in Ramsey theory.  Ron Graham, who
had a long-standing offer of 1000 USD
    for a proof or disproof of $w(k, k) < 2^{2^{2^{.^{.^{.^{2}}}}}}$, a tower of $k$ 2s, paid S. Shelah 500 USD for Shelah's
     improvement [17] of the bound obtainable from van der Waerden's
original proof, and paid T. Gowers 1000 USD for
      Gowers's upper bound.  Graham currently offers 1000 USD [3] for a proof or
disproof of $w(k, k) < 2^{k^2}$.

    Recently, there have been two breakthroughs in the study of the van der Waerden function $w(k, m)$.
    The first was the elegant proof by Graham [7] that if one defines 
$w_1(3,s)$ to be the least $n$ such that
    every $2$-coloring of $[1, n]$ gives either 
 a 3-term
     arithmetic progression in the first color or
$s$ consecutive numbers in the second color, then
\[s^{c\log s} < w_1(3,s) < s^{ds^2},\]
for suitable constants $c,d > 0$.  Of course this
immediately gives $w(k, 3) < k^{dk^2}$ since we
trivially have $w(k,3)=w(3,k) \leq w_1(3,k)$.  The second was the
amazing (computer) calculation $w(6, 6) = 1132$ by Kouril [11], extending the list of
previously known values $w(3, 3) = 9$, $w(4, 4) = 35$, and $w(5, 5)
= 178.$ A list of other known exact values of $w(k,m)$ appears in
[13]. In view of Graham's bounds on $w_1(3,s)$, it would be extremely desirable
    to obtain improved bounds on $w(k, 3)$.
      Of particular interest is the question of whether or not there is a non-polynomial lower bound for $w(k, 3)$.

    In this note we give a lower bound of $k^{(2 - o(1))} < w(k, 3)$. 
This seems weak, although we have $w(k, 3) < k^2$,
    for $5 \leq k \leq 16$  (see Table 1).
    It is our hope that others will find the question of obtaining improved bounds for $w(k, 3)$ to be interesting.
     Ultimately, one would like to find the true order of growth of the functions $w(k, 3)$, $w(k,
4), \dots, w(k, k).$ Perhaps this will be accomplished in this century, perhaps not!
Our quest for a lower bound  for $w(k,3)$ turned
(quite naturally) into a lower bound for $w(k,m)$ for an arbitrary fixed
$m$.  We also present an upper bound for $w(k,4)$, an upper bound for
$w(4;s)$, and a lower bound for $w(k;s)$ for an arbitrary fixed $k$.  (The function
$w(k;s)$ is defined below, in Section 2).
Section 2 contains the just-mentioned bounds.  In Section 3 we define several other
related functions and discuss some relationships
among these various functions.  We also provide  a table of values
of these functions for small values of  $s$ when $k=3$.

Note that we use $c$ and $d$ repeatedly to stand for positive constants,
but that these constants generally differ from paragraph to paragraph.  The context
will always make clear the meaning of a particular constant.

\section{Upper and Lower Bounds for Certain van der Waerden Functions}

We shall need several definitions, which we collect here.

For positive integers $k$ and $n$,
$$
r_k(n) = \max_{S \subseteq [1,n]} \{|S|: S \mbox{ contains no $k$-term
arithmetic progression}\}.
$$

For positive integers $k$ and $m$, denote by $\chi_k(m)$ the minimum number
of colors required to color $[1,m]$ so that there is no monochromatic
$k$-term arithmetic progression.

The function $w_1(3,s)$ has been defined in Section 1.  Similarly, we
define $w_1(k,s)$ to be the least $n$ such that every $2$-coloring of $[1,n]$ admits
either a $k$-term arithmetic progression of the first color or
$s$ consecutive integers of the second color.

Lastly, for positive integers $k$ and $s$, we denote the least positive integer $n$
such that every $s$-coloring of $[1,n]$ admits a monochromatic $k$-term
arithmetic progression by $w(k;s)$.

We begin with an upper bound for $w_1(4,s)$.  The proof is essentially the same
as the proof given by Graham [7] of an upper bound for $w_1(3,s)$.  For completeness,
we include the proof here.  We will make use of a recent result
of Green and Tao [9], who showed that for some constant $c>0$,
\begin{equation}
r_4(n) < n e^{-c \sqrt{\log \log n}}
\end{equation}
for all $n \geq 3$.

\begin{prop}  There exists a constant $c>0$ such that
$w_1(4,s) < e^{s^{c \log s}}$ for all $s \geq 2$.
\end{prop}

\noindent
{\it Proof.}  Suppose we have a $2$-coloring of
$[1,n]$ (assume $n \geq 4$) with no $4$-term arithmetic progression of the
first color and no $s$ consecutive integers of the second
color.  Let $t_1<t_2<\cdots<t_m$ be the integers
of the first color.  Hence, $m < r_4(n)$.  
Let us define $t_0=0$ and $t_{m+1}=n$.
Then there must be some $i$, $1 \leq i \leq m$, such that
$$
t_{i+1}-t_i > \frac{n}{2r_4(n)}. 
$$
(Otherwise, using $r_4(n) \geq 3$, we would have
$n = \sum_{i=0}^n (t_{i+1}-t_i) \leq \frac{n(m+1)}{2r_4(n)}
\leq \frac{n(r_4(n)+1)}{2r_4(n)} \leq \frac{n}{2} + \frac{n}{6}$.)

Using (1), we now have an $i$ with
$$
t_{i+1}-t_i > \frac{n}{2r_4(n)}>\frac{1}{2}e^{c\sqrt{\log \log n}}. 
$$

If $n \geq e^{s^{d\log s}}$, $d = c^{-2}$, then
$\frac{1}{2}e^{c \sqrt{\log \log n}} \geq s$ and we have
$s$ consecutive integers of the second color,
a contradiction.  Hence, $n < e^{s^{d \log s}}$
and we are done.
\hfill $\Box$

Clearly $w(4,s) \leq w_1(4,s)$.  Consequently, we have the
following result.

\begin{cor}
There exists a constant $d>0$ such that
$w(k,4) < e^{k^{d \log k}}$ for all $k \geq 2$.
\end{cor}

Using Green and Tao's result, it is not too difficult
to obtain an upper bound for $w(4;s)$.

\begin{prop} There exists a constant $d>0$ such that
$w(4;s) < e^{s^{d \log s}}$ for all $s \geq 2$.
\end{prop}

\noindent
{\it Proof.} Consider a $\chi_4(m)$-coloring of $[1,m]$
for which there is no monochromatic $4$-term arithmetic progression.
Some color must be used at least $\frac{m}{\chi_4(m)}$ times,
and hence $\frac{m}{\chi_4(m)} \leq r_4(m)$ so that
$\frac{m}{r_4(m)} \leq \chi_4(m)$.  Let $c>0$ be such that
(1) holds for all $n \geq 3$,
and let $m = e^{s^{d \log s}}$, where
$d = c^{-2}$.  Then $\chi_4(m) \geq \frac{m}{r_4(m)}
> e^{c \sqrt{\log \log m}}=s$.  This means that every $s$-coloring
of $[1,m]$ admits a monochromatic $4$-term arithmetic progression.
Since $m=e^{s^{d \log s}}$, the proof is complete.
\hfill $\Box$

It is interesting that the bounds in Corollary 2.2 and Proposition 2.3
have the same form.

The following theorem is deduced without too much
difficulty from
the Symmetric Hypergraph Theorem as it appears
in [8], combined with an old result of Rankin [16]. To the best of our knowledge it
has not appeared in print before,
even though
it is better, for large $s$,
than the standard bound
$\frac{cs^k}{k} (1 + o(1))$ (see [8]), the bound $s^{k+1} - \sqrt{c (k+1) \log (k+1)}$
by Erd\H{os} and Rado [4], and the bound
$\frac{ks^k}{e(k+1)^2}$ by Everts [5].  We give the proof in some detail.
The proof makes use of the following facts:
\begin{equation}
\chi_k(n) < \frac{2n \log n}{r_k(n)} (1+o(1)),
\end{equation}
which appears in [8] as a consequence of the Symmetric Hypergraph
Theorem,
and
\begin{equation}
r_k(n) > ne^{-c (\log n)^{\frac{1}{z+1}}},
\end{equation}
which, for some constant $c>0$, holds for all $n \geq 3$ (this appears in [16]).

\begin{thm}  Let $k \geq 3$ be fixed, and
let $z = \lfloor \log_2 k\rfloor$.  There exists a
constant $d>0$ such that $w(k;s) > s^{d (\log s)^z}$ for all
sufficiently large $s$. 
\end{thm}

\noindent
{\it Proof.} Fix $k \geq 3$ and let $z = \lfloor \log_2 k\rfloor$. 
Note that for positive integers $s$ and $m$,
$$
[s \geq \chi_k(m)] \Rightarrow [w(k;s) >m].
$$
This observation, which can be verified by unraveling the definitions, is
an essential ingredient of the proof.

For large enough $m$, (2) gives
\begin{equation}
\chi_k(m) < \frac{2m\log m}{r_k(m)} \left(1+\frac12 \right) = \frac{3m\log m}{r_k(m)}.
\end{equation}

Now let $d =\left( \frac{1}{2c}\right)^{z+1}$, and let $m = s^{d (\log s)^z}$, where
$s$ is large enough so that (1) holds.  By (3),
noting that
$\log m =d(\log s)^{z+1}= \left(\frac{\log s}{2c}\right)^{z+1}$, we have
$$
\frac{m}{r_k(m)} < e^{c (\log n)^{\frac{1}{z+1}}} = e^{c \cdot \frac{\log s}{2c}}
= \sqrt{s}.
$$
Therefore,
$$
\frac{3m\log m}{r_k(m)} < 3 d\sqrt{s}  (\log s)^{z+1} <s
$$
for sufficiently large $s$.  Thus, for sufficiently large $s$,
$$
\chi_k(m) < \frac{3m\log m}{r_k(m)} <s.
$$
According to the observation at the beginning of the proof, this
implies that $w(k;s) > m = s^{d (\log s)^z}$, as required.
\hfill $\Box$

We now give a lower bound on $w(k,m)$. We make use of the L\a{o}vasz
Local Lemma (see [8] for a proof), which will be implicitly stated
in the proof.

\begin{thm} Let $m \geq 3$ be fixed.  Then for all sufficiently
large  $k$,
\[ w(k,m) > k^{m-1-\frac{1}{\log \log k}} . \]
\end{thm}

\noindent
{\it Proof.}  Given $m$, choose $k > m$ large enough so that
    \begin{equation}
    k^{\frac{1}{2 m\log \log k}} > \left(m - \frac{1}{2 \log \log
k}\right)\log k
    \end{equation}
    and
    \begin{equation}
    6 < \frac{\log k} {\log \log k}.
    \end{equation}

    Next, let $n = \lfloor k^{m-1-\frac{1}{\log \log k}} \rfloor$. To
prove the theorem, we will show that
    there exists a (red, blue)-coloring of $[1,n]$  for which there is no red $k$-term
arithmetic
    progression and no blue  $m$-term arithmetic progression.

    For the purpose of using the L\a{o}vasz Local Lemma, randomly color $[1,n]$ in the following way.
     For each $i \in [1,n]$, color $i$ red with probability $p = 1-k^{\alpha-1}$ where
     \[ \alpha \stackrel{\mbox{\tiny def}}{=} \frac{1}{2m \log \log k},\]
     and color it blue with probability $1-p$.

    Let ${\cal P}$ be any  $k$-term arithmetic progression.  Then, since $1+x \leq
e^{x}$,  the probability that ${\cal P}$
     is red is
 \[ p^k = \left(1-k^{\alpha -1}\right)^{k} \leq (e^{-k^{\alpha -1}})^k = e^{-k^{\alpha}}. \]
Hence, applying (5), we have
\[ p^k < \left(\frac{1}{e}\right)^{\left(m- \frac{1}{2
\log \log k}\right)\log k} = \frac{1}{k^{m-\frac{1}{2\log \log k}}} .\]

Also, for ${\cal Q}$ any $m$-term arithmetic progression, the
probability that ${\cal Q}$ is blue is
\[ (1-p)^m = (k^{\alpha -1})^m = \frac{1}{k^{m-\frac{1}{2 \log \log
k}}} .\]
    Now let ${\cal P}_1, {\cal P}_2,\ldots, {\cal P}_t$ be  all of the arithmetic
progressions in $[1,n]$ with length
    $k$ or $m$.  So that we may apply the L\a{o}vasz Local Lemma,
      we form the ``dependency graph'' $G$  by setting $V(G)= \{{\cal P}_1,{\cal
P}_2,\ldots, {\cal P}_t\}$ and
      $E(G) = \{\{{\cal P}_i,{\cal P}_j\}: i \neq j, {\cal P}_i \cap {\cal P}_j \neq
\emptyset \}$. For each
${\cal P}_i \in V(G)$, let $d({\cal P}_i)$ denote the degree of the
vertex ${\cal P}_i$ in $G$, i.e.,  $|\{e \in E(G): {\cal P}_i \in
e\}|$.
    We now estimate $d(P_i)$ from above. Let $x \in [1,n].$ The number of  $k$-term arithmetic progressions ${\cal P}$
     in $[1,n]$ that contain $x$ is bounded above by $k\cdot\frac{n}{k-1}$, since there are $k$ positions that $x$ may occupy
     in ${\cal P}$ and since the gap size of ${\cal P}$ cannot exceed
$\frac{n}{k-1}$. Similarly, the number
     of $m$-term arithmetic progressions ${\cal Q}$ in $[1,n]$ that contain $x$ is bounded above by $m \cdot \frac{n}{m-1}$.

    Let ${\cal P}_i$ be any $k$-term arithmetic progression contained in $[1,n]$.  The total number of
     $k$-term arithmetic progressions ${\cal P}$ and  $m$-term arithmetic progressions ${\cal Q}$ in $[1,n]$
       that can intersect ${\cal P}_i$  is bounded above by
       \[k\left(k \cdot \frac{n}{k-1}+ m \cdot \frac{n}{m-1}\right)  <
kn\left(2+\frac{2}{m-1}\right), \]
         since $k > m$.  Thus, $d({\cal P}_i)   < 
kn\left(2+\frac{2}{m-1}\right)$ when $|{\cal P}_i|=k$.
    Likewise, $d({\cal P}_i) < mn\left(2+\frac{2}{m-1}\right)$ when
$|{\cal P}_i| = m$.
    Thus, for all vertices ${\cal P}_i$   of $G$, we have $d({\cal P}_i) <
kn\left(2+\frac{2}{m-1}\right)$.

    To finish setting up the hypotheses for the L\a{o}vasz Local Lemma,
we let $X_i$   denote the event that the arithmetic
     progression ${\cal P}_i$   is
 \[ \left\{ \begin{array}{ll}
 \mbox{ red } & \mbox{ if $|{\cal P}_i| =k$} \\
 \mbox{ blue } & \mbox{ if $|{\cal P}_i| = m$.}
 \end{array} \right. \]

We have seen above that for all $i$, $1 \leq i \leq t$, the
probability of the event $X_i$ is less than
\[ q \stackrel{\mbox{\tiny def}}{=} \frac{1}{k^{m - \frac{1}{2 \log \log
k}}}.
\]

Let $d = \displaystyle \max_{1 \leq i \leq t} d({\cal P}_i).$ We
showed above that
\[ d < 2kn\left(1+\frac{1}{m-1}\right) .\]

    We are now ready to apply the L\a{o}vasz  Local Lemma, which says
that in these circumstances,
    if the condition $eq(d+1) < 1$ is satisfied, then there is a (red, blue)-coloring
of $[1,n]$
     such that no event
    $X_i$  occurs, i.e., there is a (red, blue)-coloring of $[1,n]$  for which there
is no red  $k$-term arithmetic
    progression and no blue  $m$-term arithmetic progression.  This will imply
    \[ w(k,m) > n  = k^{m-1-\frac{1}{\log \log k}}, \]
      as desired.
    Thus, the proof will be complete when we verify that $eq(d+1) < 1$.
    Using  $m\geq 3$, we have $d < 3kn$, so that $d+1 < 3kn+1 < e^{2}kn$.  Hence, it is sufficient to verify that
    \begin{equation}
    e^{3}qkn < 1.
\end{equation}
      Since $q= \frac{1}{k^{m-\frac{1}{2 \log \log k}}}$ and $n \leq k^{m-1-\frac{1}{\log
\log k}}$,
    inequality (7) may be reduced to (6), and the proof is now complete.
\hfill $\Box$

 \noindent {\it Remark.} For condition (5), it suffices to have
$k^{\frac{1}{2 m\log \log k}} > m \log k$, or $\log k > 2m \log m
(\log \log k) + 2m(\log \log k)^{2}$.  
When $k \geq e^{e^{m^3}}$, this condition becomes
$e^{m^3} > 2m^4 \log m + 2m^7$.  Since, for $m \geq 3$,
we have $e^{m^3} > m^9 > 2m^4 \log m + 2m^7$, having $k >
e^{e^{m^3}}$ is sufficient for both (5) and (6).

\section{Some Related Functions}
In this section we define some functions related to $w(k,m)$  and
mention various bounds for, and relationships among, these.
For reference, we define all functions used in this section
(including those already defined).


\noindent ${\bf w(k,s)}$ is the least positive integer $n$
such that every $2$-coloring of $[1,n]$ admits either 
a $k$-term
arithmetic progression of the first color
or
an $s$-term
arithmetic progression of the second color.

\noindent ${\bf w_1(k,s)}$ is the least positive integer $n$
such that every $2$-coloring of $[1,n]$ admits either 
 a $k$-term
arithmetic progression of the first color or 
$s$ consecutive integers of the second color.

 \noindent ${\bf w(k;s)}$ is the least positive integer $n$
such that every $s$-coloring of $[1,n]$ admits a monochromatic
$k$-term arithmetic progression. 

 \noindent 
 ${\bf G(k,s)}$ is the least positive
integer $n$ such that for every set $S=\{x_1 < x_2 < \cdots < x_n\}$
with $x_i - x_{i-1} \leq s$, $2 \leq i \leq n$, $S$ contains a
$k$-term arithmetic progression. 

 \noindent 
 ${\bf M(k,s)}$ denotes the least
positive
integer $n$ such that whenever $X= \{x_1,x_2,\ldots,x_n\}$ and 
$x_i \in [(i-1)s, is-1]$, $1 \leq i \leq n$, there is a $k$-term arithmetic progression
in $X$.

 \noindent 
${\bf w^{*}(k;s)}$ denotes the least positive integer $n$ such that
every $s$-coloring  $\chi: [1,n] \rightarrow [1,s]$ admits either a
monochromatic $k$-term arithmetic progression or a $k$-term
arithmetic progression whose colors form an arithmetic progression
(increasing or decreasing).


We start with the following inequalities
involving $w_1(k,s)$.

\begin{prop} For any positive integers $k$ and $s$,
the following hold:
\begin{enumerate}
\item[(i)] $w_1(k,s) \leq s M(k,s)$;
\item[(ii)] $w_1(k,s) \leq sG(k,s)$;
\item[(iii)] $w_1(k,s) \leq w(k;s)+s$.
\end{enumerate}
\end{prop}
 {\em Proof.}
As the proofs of (i) and (ii) are quite similar, we include
the proof of (i) and leave the other to the
 reader. Let $m = M(k,s)$ and let $n=sm$. Let $\chi$ be any
(red, blue)-coloring of $[1,n]$. Assume there are no $s$ consecutive
blue integers. So, for each $i$,  $1 \leq i \leq m$, the interval
$[(i-1)s+1,is]$ contains a red element, say $a_i$. Then, by the
definition of $M(k,s)$, there is a $k$-term arithmetic progression among
the
$a_i$'s.

We now show (iii).
By definition, there exists a (red, blue)-coloring of $[1,w_1(k,s)-1]$
with no red $k$-term arithmetic progression and no $s$ consecutive blue elements.
Let the red elements under this coloring be
$\mathcal{R} = \{r_1<r_2<\dots<r_t\}$.  Note that
$r_1 \leq s$ and $r_t \geq w_1(k,s)-s-1$.  Define the
following $s$-coloring of
$[r_1,w_1(k,s)-1]$.
Color all elements in $\mathcal{R}$ with color $0$.
For $i=1,2,\dots,s-1$, in order,
color all elements in $(\mathcal{R}+i) \setminus \bigcup_{j=0}^{i-1}
(\mathcal{R}+j)$ with color $i$.  This is well defined since
$r_{x+1}-r_x \leq s$ for any $x$.
Since $\mathcal{R}$ contains no $k$-term arithmetic progression, none of
$(\mathcal{R}+i) \setminus \bigcup_{j=0}^{i-1}
(\mathcal{R}+j)$ contain a $k$-term arithmetic progression.  
Since $[r_1,w_1(k,s)-1]$ contains at least
$w_1(k,s)-s-1$ elements and our
$s$-coloring admits no monochromatic $k$-term arithmetic progression,
we see that
 $w(k;s) \geq w_1(k,s) - s$.
\hfill $\Box$

\noindent
{\it Remark.} Using $w(4,k) \leq w_1(4,k)$ and part (iii) from the
above proposition, we see that Proposition 2.1 and
Corollary 2.2 follow from
Proposition 2.3, without appealing to Graham's argument.

In the next proposition, we give an alternate way of describing
$M(3,s)$. Before doing so, we introduce some terminology. A {\em
3-term ap$^{+}$} is an ordered triple of the form $x,x+d,x+2d+1$
where $x,d \geq 1$. A {\em 3-term $ap^{-}$} is an ordered triple of
the form $x,x+d,x+2d-1$, where $x \geq 1$ and $d \geq 2$. Note that,
in either case, the three terms are distinct. Finally, an ordered triple
$a,b,c$ is called an {\it arithmetic progression(mod $s$)} if $c-b \equiv b-a$ (mod
$s$).

We make use of the following lemma.

\begin{lemma}
Let $s \geq 2$. For $i\in \mathbb{Z}^+$, let $B_{i} = [(i-1)s+1,is]$. For
$j\in \mathbb{Z}^+$, define $r_j$ to be the unique member of $\{1,2,\dots,s\}$
such that $j \equiv r_j$ (mod $s$). Let $1 \leq x < y < z$ with $x
\in B_{i_1}$, $y \in B_{i_2}$, and $z \in B_{i_3}$. Then $x,y,z$ is
an arithmetic progression if and only if one of the following holds:
\begin{enumerate}
\item[(i)] $r_z - r_y = r_y-r_x=a$  and
 $i_3-i_2=i_2-i_1=b$ where either (i) $a > 0$, or (ii) $a \leq 0$ and $b > 0$;
 \item[(ii)] 
\begin{enumerate} 
\item  $r_x,r_y,r_z$, with $r_x<r_y$,
 is an arithmetic progression (mod $s$), but not an arithmetic progression, 
and
\item
 either $i_1,i_2,i_3$ is a $3$-term $ap^{+}$ or $i_1=i_2 = i_3-1$;
\end{enumerate}
 \item[(iii)] \begin{enumerate} \item
 $r_x,r_y,r_z$, with $r_x>r_y$, is an arithmetic progression (mod $s$), but not an
arithmetic progression, and 
\item either $i_1,i_2,i_3$ is a $3$-term $ap^{-}$ or
 $i_3=i_2=i_1+1$. 
\end{enumerate} 
\end{enumerate}
 \end{lemma}

\noindent {\em Proof.} If (i) holds, then for some $d$ satisfying
$|d| \leq \lfloor \frac{s-1}{2}\rfloor$ we have \[ z-y =
(i_3-i_2)s+d = (i_2 -i_1)s+d = y-x.\] 

Now assume (ii) holds. Note
that $r_y-r_x = s + r_z - r_y.$ Thus, \noindent $z-y = (i_3-i_2)s -
r_y + r_z = (i_2-i_1+1)s - r_y+r_z = (i_2-i_1)s + (s-r_y+r_z) =
(i_2-i_1)s + r_y - r_x = y - x$.

Now assume (iii) holds. In this case $r_z-r_y = s+r_y-r_x$. Therefore,
$$z-y = (i_3-i_2)s + r_z - r_y = (i_2 - i_1)s - s +r_z - r_y = (i_2 - i_1)s + r_y -
r_x = y -x.$$

For the converse, it suffices to consider three cases. 

\noindent
{\tt Case 1.} $r_x \leq r_y \leq r_z$ or $r_x \geq r_y \geq r_z$. 
In this case, it is clear that $i_3-i_2 = i_2 - i_1$.

\noindent 
{\tt Case 2.} $r_x < r_y$ and $r_z \leq r_y$. In this case, $i_3 -
i_2 > i_2 -i_1$. Furthermore, $i_3-i_2 \geq i_2 - i_1 +2$ is not
possible, since then we would have $z-y \geq y-x + s+1$. Hence
$i_3 -i_2 = i_2 - i_1 + 1,$
so that $i_1,i_2,i_3$ is an $ap^{+}$. 
Also, \[ z-y = (i_3-i_2)s - (s-r_y+r_z),\] and \[ y-x = (i_2-i_1)s -
(r_y-r_x).\] So $r_y-r_x \equiv r_z-r_y$(mod $s$).

\noindent 
{\tt Case 3.} $r_x > r_y$ and $r_z \geq r_y$. The proof is almost
the same as that for Case 2, and we leave it to the reader.
 \hfill $\Box$

\begin{prop}  For all $s\geq 2$, $M(3,s)$ is the least positive integer $n$
such that every $s$-coloring $\chi: [1,n] \rightarrow [1,s]$ admits
a triple $A = \{a < b < c\}$ satisfying one of the following:
\begin{enumerate} 
\item[(i)] $A$ is an arithmetic
progression and $\chi(b)-\chi(a) = \chi(c) - \chi(b)$ (possibly
negative or 0);
\item[(ii)] $A$ is a $3$-term $ap^{+}$ and
$(\chi(a),\chi(b),\chi(c))$, with $\chi(a) < \chi(b)$, is an arithmetic progression
(mod
$s$), but not an arithmetic progression;
\item[(iii)] $A$ is a $3$-term $ap^{-}$ and
$(\chi(a),\chi(b),\chi(c))$, with $\chi(a) > \chi(b)$, is an arithmetic progression
(mod $s$), but not a (decreasing) arithmetic progression.
\end{enumerate} 
\end{prop} 

\noindent {\em Proof.} Let an $s$-coloring of $[1, n]$ be given, using the colors
$1,2,\dots, s$. We use this coloring to define a set $\{x_1, x_2, \dots , x_n\}$ as
follows:  For each $i$, let $x_i$ be that element of the block $B_i =
[(i-1)s+1, is]$ which is congruent to the color of $i$.  By Lemma
3.2, the minimum $n$ such that any set $\{x_1, x_2, \dots, x_n\}$
constructed in this way must contain a 3-term arithmetic progression is $M(3,s)$.
\hfill $\Box$

The following  inequalities are proved by Nathanson [14].

\begin{thm}(Nathanson) For all positive integers $k$ and $s$,
\begin{enumerate}
\item[(i)] $G(k,s) \leq s M(k,s)$,

\item[(ii)] $M(k,s) \leq G(k,2s-1)$,

\item[(iii)] $G(k,s) \leq w(k;s)$,

\item[(iv)] $M(k,s) \leq w(k;s)$,

\item[(v)] $w(k;s) \leq M(s(k-1)+1,s)$, and

\item[(vi)] $w(k;s) \leq G(s(k-1)+1,2s-1)$.

\end{enumerate}
\end{thm}

Investigating $G(k,s)$, Alon and Zaks [1] have shown the following
result.  Evaluating their bound when $s=2$ gives, to our eyes, 
a surprising result since we may view this (loosely speaking) as a
$2$-coloring of $[1,2^{k(1+o(1))}]$ with no monochromatic $k$-term
arithmetic progression, where one of the color classes has no two
consecutive integers.  It is surprising that they proved a
lower bound for such a resticted family of colorings that is almost
as large as the best known lower bound for $w(k,k)$.

\begin{thm} (Alon and Zaks)  For every $s \geq 2$, there exists
a constant $c>0$ (dependent upon $s$) such that
$$
G(k,s)>s^{k - c \sqrt{k}}
$$
for all $k \geq 3$.
\end{thm}

Note that, from part (vi) of Theorem 3.4, an upper bound for
$G(k,s)$ would give an upper bound for $w(k;s)$.  In particular,
an upper bound on $G(k,3)$ for odd $k$ would lead to an upper bound for
$w(k,k)$.

We now offer a few more inequalities involving the
functions discussed in this section.  We use
Szemeredi's result on arithmetic progressions [18]
to prove part (iii) of Proposition 3.6, which improves
(for large $s$) part (i) of Theorem 3.4.

\begin{prop}
 For all positive integers $k$ and $s$,
the following hold.
\begin{enumerate}
\item[(i)]
 $w(k,s) \leq w_1(k,s)$
\item[(ii)]
$M(k,s) \leq w^{*}(k;s) \leq w(k;s)$
 \item[(iii)] Let $c>0$ be constant.  For $k$ fixed and $s$ sufficiently large,
$G(k,s) < csM(k,s)$.
\item[(iv)] $w_1(k,2s-1) \geq s (M(k,s)-1)+1$

\end{enumerate}
\end{prop}
 {\em Proof}. Parts (i) and (ii) are immediate from the definitions and Proposition
3.3.
 To prove (iii), we assume
that
$cs$ is an integer and
let $t = csM(k,s)$. Let
$X=\{x_1<x_2< \cdots<x_t\}$ be a set such
that $x_{i}-x_{i-1} \leq s$ for $i=2,3,\dots,t$.
We may assume $x_1 \in [1,s]$.
We must show that $X$ contains a $k$-term arithmetic progression.
Define  $X_i=\{x_{(i-1)cs+1}, \dots, x_{ics}\}$
for $i=1,2,\dots, M(k,s)$.  If for some $i$ we have
$X_i \subseteq [(i-1)s+1,is]$ then,
for
$s$  sufficiently large,
Szemeredi's result on arithmetic progressions tells us
that $X_i$ must contain a $k$-term arithmetic progression.
Assuming this is not the case, we must then have that each interval $[(i-1)s+1,is]$
for $1 \leq i \leq M(k,s)$ contains at least one element of $X$.
By the definition of $M(k,s)$, we have our $k$-term arithmetic
progression.

We now prove (iv).
We will show that if $s$ is even, then
\[ w_1(k,s-1) \geq \frac{s}{2} \left(M\left(k, \frac{s}{2}\right) -1 \right) + 1.
\]
The case when $s$ is odd is similar and is left to the reader.
 Let $n = M\left(k,\frac{s}{2}\right)$. Then there exists
$X=\{x_1,x_2,\ldots,x_{n-1}\}$ containing no $k$-term arithmetic
progression and such that $x_i \in \left[\frac{(i-1)s}{2}+1,\frac{is}{2}\right]$ for $1
\leq i
\leq n-1$. Consider the following 2-coloring of $\left[1,\frac{(n-1)s}{2}\right]$:
$\chi(x)=0$ for each $x_i \in X$, and $\chi(x)=1$ otherwise. Clearly
there is no $k$-term arithmetic progression with color 0. Also, since
each interval $\left[\frac{(i-1)s}{2}+1,\frac{is}{2}\right]$ contains an element with
color 0, the longest string of consecutive elements with color 1 has length
not exceeding $s-2$, which implies the desired result. \hfill $\Box$

We end with a table of computed values.  These were all computed
with a standard backtrack algorithm except for
$w(3,14), w(3,15),$ and $w(3,16)$, which are due to
Michal Kouril [12].  Based on the values in this table,
we make the following conjecture (only the last inequality
is known to hold).

\noindent
{\bf Conjecture}  For all $s \geq 2$,
$$
G(3,s) \leq w(3,s) \leq M(3,s) \leq w_1(3,s) \leq w^*(3,s) \leq w(3;s).
$$

$$
\begin{array}{r||c|c|c|c|c|c|c|c|c|c|c|c|c|c|c|}
s &2&3&4&5&6&7&8&9&10&11&12&13&14&15&16\\ \hline
G(3,s)&5&9&11&17&22&33&37&48&?&?&?&?&?&?&?\\ \hline
w(3,s)&6&9&18&22&32&46&58&77&97&114&135&160&186&218&238\\ \hline
M(3,s)&7&11&18&29&37&48&?&?&?&?&?&?&?&?&?\\ \hline
w_1(3,s)&9&23&34&73&113&193&?&?&?&?&?&?&?&?&?\\ \hline
w^*(3,s)&9&23&40&\geq 75&?&?&?&?&?&?&?&?&?&?&?\\ \hline
w(3;s)&9&27&76&?&?&?&?&?&?&?&?&?&?&?&? \\ \hline
\end{array}
$$
\centerline{\bf Table 1:  Small values of van der Waerden-like functions}

\vskip 20pt
\noindent{\bf References}

\parindent=0pt \footnotesize
\begin{enumerate}

\item N. Alon and A. Zaks, Progressions in sequences of nearly consecutive integers,
{\it J. Comb. Theory, Series A} {\bf 84} (1998), 99-109.

\item E. Berlekamp, A construction for partitions which avoid long
arithmetic progressions, {\it  Canad. Math. Bull.} {\bf  11} (1968), 409-414.

\item F. Chung, P. Erd\H{o}s, and R. Graham, On sparse sets hitting
linear forms, Number theory for the millennium, I (Urbana, IL,
2000), 257-272, A K Peters, Natick, MA, 2002.

\item P. Erd\H{o}s and R. Rado, Combinatorial theorems on classifications of
subsets of a given set, {\it Proc. London Math. Soc.} {\bf 3} (1952), 417-439.

\item F. Everts, Colorings of sets, Ph.D. thesis, University of Colorado, 1977.

\item W. T. Gowers, A new proof of Szemeredi's theorem, {\it Geom. Funct.
Anal.} {\bf 11} (2001), no. 3, 465-588.

\item  R. Graham,  On the growth of a van der Waerden-like function,
{\it Integers: El. J. Combinatorial
Number Theory} {\bf 6} (2006), A29.

\item   R. Graham, B. Rothschild, and
J. Spencer, \underline{Ramsey Theory}, Wiley-Interscience,
2e, 1990.

\item B. Green and T. Tao, New bounds for Szemeredi's theorem II: A new bound for
$r_4(N)$, preprint:  {\tt arXiv:math//0610604v1}.

\item J. H. Kim, The Ramsey number
$R(3,t)$ has order of magnitude $t\sp 2/\log t$, {\it Random Structures
Algorithms} {\bf 7} (1995), no. 3, 173-207.

\item M. Kouril, Ph.D. thesis, University of Cinncinati, 2007.

\item M. Kouril, private communication, 2007.

\item B. Landman, A. Robertson, and C. Culver,
Some new exact van der Waerden numbers,
{\it Integers:  El. J. Combinatorial Number Theory} {\bf 5(2)}
(2005), A10.

\item M. Nathanson, Arithmetic progressions contained in sequences
with bounded gaps, {\it Canad. Math. Bull.} {\bf 23} (1980), 491-493.

\item F. Ramsey, On a problem of formal logic, {\it Proc. London Math.
Soc.} {\bf 30} (1930), 264-286.

\item R. Rankin, Sets of integers containing not more than a given numbers
of terms in arithmetical progression, {\it Proc. Roy. Soc.
Edinburgh Sect. A} {\bf 65} (1960/1961), 332-344.

\item S. Shelah, Primitive recursive bounds for van der Waerden numbers,
{\it J. American Math. Soc.} {\bf 1} (1988), 683-697.
\item E. Szemeredi, On sets of integers containing no $k$ elements in
arithmetic progression, {\it Acta Math. Acad. Sci. Hungaricae} {\bf 20} (1975),
199-245.

\item B. L. Van der Waerden, Beweis einer baudetschen Vermutung,
{\it Nieuw Archief voor Wiskunde} {\bf 15} (1927), 212-216.

\end{enumerate}

 \end{document}